\documentclass[11pt,a4paper]{article}
\usepackage{amsmath,amssymb}
\newtheorem{theorem}{Theorem}[section]

\newtheorem{corollary}{Corollary}[section]
\newtheorem{definition}{Definition}[section]
\newtheorem{proposition}{Proposition}[section]
\newenvironment{remark}{\vspace{1ex}{\bf Remark.}\rm}{\vspace{1ex}}
\newenvironment{proof}{{\bf Proof.}}{\par\hspace{25em}\rule{1ex}{1ex}\par}
\newcommand{\ds}{\displaystyle}

\title{Full description of totally geodesic unit vector fields\\
on 2-dimensional Riemannian manifolds.}
\author{Yampolsky A.}
\date{}
\begin{document}
\maketitle
\begin{abstract}
    We give a full geometrical description of local totally geodesic unit vector field on Riemannian 2-manifold,
    considering the field as a local imbedding of the manifold into its unit tangent bundle with the Sasaki metric.

    {\it AMS subject class:} Primary 54C40,14E20; Secondary 46E25, 20C20
\end{abstract}

\section*{Introduction}
Let $(M,g)$ be an $(n+1)$ -- dimensional Riemannian manifold  with metric $g$. A vector field $\xi$ on it is
called {\it holonomic} if $\xi$ is a field of normals of some family of regular hypersurfaces in $M$ and {\it
non-holonomic} otherwise. The foundation of the classical geometry of unit vector fields was proposed by A.Voss
at the end of the nineteenth century. The theory includes the {\it Gaussian} and {\it the mean curvature} of a
vector field and their generalizations (see \cite{Am} for details).

Recently, the geometry of vector fields has been considered from another point of view. Let $T_1M$ be a unit
tangent bundle of $M$ endowed with the Sasaki metric \cite{S}. If $\xi$ is a unit vector field on $M$, then one
may consider $\xi$ as a mapping $\xi : M \to T_1M $. The image $\xi (M)$ is a submanifold in $T_1M$ with metric
induced from $T_1M$ and one may apply the methods from the study of the geometry of submanifolds to determine
geometrical characteristics of a unit vector field. A unit vector field $\xi$ is said to be {\it minimal}  if $
\xi(M)$ is a minimal submanifold in $T_1M$. A unit vector field on $S^3$ tangent to the fibers of the Hopf
fibration $S^3 \stackrel{S^1}{\longrightarrow}S^2$ is a unique unit vector field with globally minimal volume
\cite{G-Z}. This result fails in higher dimensions. A lower volume is achieved by a vector field with one
singular point, namely the inverse image under stereographic projection inverse image of a parallel vector field
on $E^n$ \cite{Ped}. The lowest volume is reached for the \textit{North-South} vector field with two singular
points \cite{B-Ch-N}.

 A local approach to minimality of unit vector fields was developed in \cite{GM-LF}. A number
of examples of locally minimal unit vector fields was found  [2~--~4, 7 -- 10, 12 -- 14,  16 -- 18] on various
manifolds. In \cite{Ym1} the author presented {\it an explicit expression} for the second fundamental form of
$\xi(M)$ and found some examples of vector fields with constant mean curvature. This expression is the key to
solving a problem about \textit{totally geodesic vector fields} on a given Riemannian manifold.  Originally, the
problem of a full description of \textit{all totally geodesic submanifolds} in the tangent (sphere) bundle of
spaces  of constant curvature was posed by A.Borisenko in \cite{B-Y}. The totally geodesic vector fields form a
special class of such submanifolds. In \cite{Ym2} this problem was solved in the case of 2-manifolds of constant
curvature. In \cite{Ym4} an example of a totally geodesic unit vector field was found on a surface of revolution
with non-constant but sign-preserving Gaussian curvature.

In this paper, we completely determine the Riemannian 2-manifolds which admit a unit vector field $\xi$ such
that $\xi(M)$ is a totally geodesic submanifold in $T_1M$. Moreover, we explicitly determine the vector field.
Under some restrictions, we find an isometric immersion of the metric into Euclidean 3-space which gives a
surface with the necessary properties.

\section{The main result}

Let $\xi$ be a unit vector field on a Riemannian manifold $(M^n,g)$. Then $\xi$ can be considered as a mapping
$\xi:M^n\to T_1M^n$. In this way one can use geometrical properties of the submanifold $\xi(M^n)$ to determine
the geometrical characteristics of the vector field.
\begin{definition} A unit vector field on Riemannian manifold $M^n$ is said to be totally geodesic, if the
submanifold $\xi(M^n)\subset T_1M^n$  is totally geodesic in the unit tangent bundle with the Sasaki metric.
\end{definition}
\begin{definition}
A point $q\in M^n$ is said to be stationary for the vector field $\xi$ if $\nabla_X\xi\,|_q=0$ for all $X\in
T_qM^n$.
\end{definition}
If stationary points fills a domain $D\subset M^n$, then locally $M^n=M^{n-k}\times E^k$, where $E^k$ is a
Euclidean factor of dimension $k\geq 1$. In the case $n=2$, the manifold is then flat in $D$. If the manifold is
of sign-preserving Gaussian curvature, then we can always restrict our considerations to the domain with no
stationary points of a given unit vector field. The main result of the paper is the following theorem.
\begin{theorem}
Let $M^2$ be a Riemannian manifold with sign-preserving Gaussian curvature $K$. Then , on some open subset U of
M, there exists a unit totally geodesic vector field $\xi$ if and only if
\begin{itemize}
\item[(a)] the metric $g$ on $U$ is locally of the form
$$
ds^2=du^2+\sin^2 \alpha(u) \, dv^2,
$$
where $\alpha(u)$ solves the differential equation $\quad\ds \frac{d\alpha}{du}=1-\frac{a+1}{\cos\alpha}; $
\item[(b)] the totally geodesic unit vector field  $\xi$ is of the form
$$
\xi=\cos (av+\omega_0)\,\partial_u+\frac{\sin (av+\omega_0)}{\sin\alpha(u)}\,\partial_v,
$$
where $a,\omega_0=const$.
\end{itemize}
\end{theorem}
\begin{remark} The Gaussian curvature $K$ of the metric is
\begin{equation}\label{curv}
K=\frac{d\alpha}{du}.
\end{equation}
Therefore, $\alpha(u)$ is the total curvature of the manifold along the meridian of the metric. The vector field
is parallel along meridians and bends along parallels with constant angle speed $a$ with respect to the
coordinate frame.
\end{remark}

\begin{proof}
Let $\xi$ be a given unit vector field on Riemannian manifold $M^n$. For dimension reasons, the kernel of the
linear operator $\nabla_X\xi:TM^n\to \xi^\perp$ is not empty. Therefore, there is a non-zero vector field $e_0$
such that $\nabla_{e_0}\,\xi=0$. In the case $n=2$, the field $e_0$ can be found explicitly. Denote by $\eta$ a
unit vector field on $M^2$ which is orthogonal to $\xi$. Set
$$
\nabla_\xi\xi=k\,\eta,\quad \nabla_\eta\eta=\varkappa\,\xi,
$$
where $k$ and $\varkappa$ are the signed geodesic curvatures of the integral trajectories of the fields $\xi$
and $\eta$ respectively. Introduce an orthonormal frame
$$
\ds e_0=\frac{\varkappa}{\lambda}\,\xi+\frac{k}{\lambda}\,\eta, \quad \ds
e_1=\frac{k}{\lambda}\,\xi-\frac{\varkappa}{\lambda}\,\eta, \quad \lambda=\sqrt{k^2+\varkappa^2}.
$$
The fields $e_0$ and $e_1$ are correctly defined on an open subset $U\subset M^2$ where the field $\xi$ has no
stationary points, i.e., points where $\lambda=0$. Restrict ourselves to this open part. It is elementary to
check that
\begin{equation}\label{Eq1}
\nabla_{e_0}\xi=0,\quad \nabla_{e_1}\xi=\lambda \eta.
\end{equation}

Denote by $\omega$ the angle function between $\xi$ and $e_0$. Then
\begin{equation}\label{lambda}
k=\lambda\sin\omega,\quad \varkappa=\lambda\cos\omega
\end{equation}
and we can set

\begin{equation}\label{fields}
\begin{array}{l}
\xi=\cos \omega \, e_0+\sin \omega\,e_1,\\[1ex]
\eta=\sin \omega \, e_0-\cos \omega\,e_1.
\end{array}
\end{equation}

Denote by $\mu$ and $\sigma$ the {\it signed} geodesic curvatures of the integral curves of the fields $e_0$ and
$e_1$ respectively. Then
$$
    \nabla_{e_0}\,e_0=\mu\,e_1,\quad \nabla_{e_1}\,e_1=\sigma\,e_0.
$$
In these terms, the second fundamental form of the submanifold $\xi(M)\subset T_1M$ can be expressed as
\cite{Ym2}
\begin{equation}\label{Forms}
   \Omega=
   \left[\begin{array}{cc}
           \ds -\mu\,\frac{\lambda}{\sqrt{1+\lambda^2}} & \ds \frac12\left(\sigma\,\lambda +\frac{1-\lambda^2}{1+\lambda^2}e_0(\lambda)\right)\\[2ex]
           \ds \frac12\left(\sigma\,\lambda +\frac{1-\lambda^2}{1+\lambda^2}e_0(\lambda)\right)& \ds e_1\left(\frac{\lambda}{\sqrt{1+\lambda^2}}\right)
   \end{array}\right].
\end{equation}
Set
$$\cos(\alpha/2)=\frac{1}{\sqrt{1+\lambda^2}}.
$$
Then we have
$$
\begin{array}{l}
\ds \frac{\lambda}{\sqrt{1+\lambda^2}}=\sin (\alpha/2),\quad
\frac{1-\lambda^2}{1+\lambda^2}=\cos\alpha, \\[2ex]
\ds e_0(\lambda)=\frac{e_0(\alpha)}{2\cos^2(\alpha/2)},\quad
e_1\left(\frac{\lambda}{\sqrt{1+\lambda^2}}\right)=\frac12\cos(\alpha/2)\, e_1(\alpha).
\end{array}
$$
After these simplifications
$$
   \Omega=
   \frac12\left[\begin{array}{cc}
           \ds -2\mu\,\sin(\alpha/2) & \ds \frac{\sigma\,\sin\alpha +e_0(\alpha)\cos\alpha}{2\cos^2(\alpha/2)}\\[2ex]
            \ds \frac{\sigma\,\sin\alpha +e_0(\alpha)\cos\alpha}{2\cos^2(\alpha/2)}& \ds \cos(\alpha/2)\, e_1(\alpha)
   \end{array}\right].
$$

Set $\Omega\equiv 0$. Then $\mu\equiv0$, since $\sin(\alpha/2)\equiv0$ implies $\lambda\equiv0$, which
contradicts the hypothesis. Therefore, a if totally geodesic vector field exists, then the \textit{integral
trajectories of the field $e_0$  are geodesics}.

Since $\cos(\alpha/2)\ne 0$, then
\begin{equation}\label{e1}
e_1(\alpha)\equiv0.
\end{equation}
 Introduce a local semi-geodesic coordinate system $(u,v)$ such that
$$
\partial_u=e_0,\ \partial_v=f(u,v)\,e_1,
$$
where $f(u,v)$ is some non-zero function. Then the line element of $M^2$ can be written as
$$
ds^2=du^2+f^2dv^2
$$
The condition \eqref{e1} implies $\partial_v\alpha=0$, which means that $\alpha=\alpha(u)$.

Consider now the last condition
$$
 \sigma\,\sin\alpha +e_0(\alpha)\cos\alpha=0.
$$
If $\cos\alpha\equiv0$, then $\sin\alpha\equiv 1$ and hence $\sigma\equiv 0$. This
means that $e_0$ is a parallel vector field on $M^2$ and hence $K=0$ again. Set
$$
\sigma\,\tan\alpha +e_0(\alpha)=0.
$$
With respect to the chosen  semi-geodesic coordinate system, $\sigma=-\partial_u f/f $ and we come to the
following relation
$$
\frac{\partial_u f}{f}=\cot\alpha\ \partial_u \alpha .
$$
Because of \eqref{e1}, we have  $\alpha=\alpha(u)$ and the equation above has an evident solution
$$
f(u,v)=C(v)\sin\alpha,
$$
where $C(v)\ne0$ is a constant of integration. Making  a $v$- parameter change one can always set $C(v)\equiv1$.
Therefore, the line element of a 2-manifold $M$ which admits a totally geodesic vector unit field is necessarily
of the form
\begin{equation}\label{met}
ds^2=du^2+\sin^2\alpha(u)\, dv^2.
\end{equation}

Turn now to the vector field. A direct computation yields
$$
\begin{array}{l}
\nabla_{e_0}\xi=\nabla_{e_0}(\cos\omega \, e_0+\sin\omega\, e_1)=(-e_0(\omega)-\mu)\,\eta, \\[1ex]
\nabla_{e_1}\xi=\nabla_{e_1}(\cos\omega \, e_0+\sin\omega\, e_1)=(-e_1(\omega)+\sigma)\,\eta.
\end{array}
$$
Since $\mu=0$ and $\nabla_{e_0}\xi=0$, we see that $\partial_u \omega=0$ and hence $\omega=\omega(v)$. The
second equality means, that
$$
-e_1(\omega)+\sigma=\tan(\alpha/2).
$$
With respect to a chosen coordinate system, we have
$$
\sigma=-\cot \alpha\ \partial_u \alpha
$$
and hence
$$
\partial_v\omega=\sin\alpha\,(\sigma-\tan (\alpha/2)=-\cos \alpha\ \partial_u\alpha-2\sin^2(\alpha/2)
$$
The right hand side does not depend on the $v$- parameter and therefore $\partial^2_{vv}\omega=0$ which means
that
$$
\omega=av+\omega_0, \quad (a,\omega_0=const).
$$
As a consequence, we come to the following differential equation for the function $\alpha(u)$:
$$
\cos \alpha\ \partial_u\alpha +2\sin^2(\alpha/2)=-a
$$
or equivalently
\begin{equation}\label{alf}
\frac{d\alpha}{du} =1-\frac{a+1}{\cos\alpha}.
\end{equation}

The proof is complete.
\end{proof}

\begin{remark}
A direct computation shows that if $\alpha$ is a solution of \eqref{alf}, then Gaussian curvature of the metric
\eqref{met} takes the form \eqref{curv}. Since it is supposed that $K$ is sign-preserving, the relation
\eqref{curv} allows to choose $\alpha$ as a new parameter on $u$-curves. With respect to the parameter $\alpha$
we have
$$
du=\frac{d\alpha}{K}=-\frac{\cos\alpha}{a+1-\cos\alpha}\,d\alpha
$$
and the line element \eqref{met} takes the form
\begin{equation}\label{met1}
ds^2=\left(\frac{\cos\alpha}{a+1-\cos\alpha}\right)^2\,d\alpha^2+\sin^2\alpha \,dv^2.
\end{equation}
\end{remark}
\begin{remark} If $\xi$ is a unit vector field on the Riemannian manifold $M^n$, then the induced metric on $\xi(M^n)$
is $ d\tilde s^2=g_{ik}du^idu^k+\big<\nabla_i\xi,\nabla_k\xi\big>du^i du^k. $ If $\xi$ is a totally geodesic
vector field on $M^2$, then the metric of $M^2$ has the standard form \eqref{met} and
$\nabla_{\partial_u}\xi=\nabla_{e_0}\xi=0,\ \nabla_{\partial_v}\xi=\sin\alpha
\nabla_{e_1}\xi=\sin\alpha\lambda\, \eta=2\sin^2(\alpha/2)\,\eta.$ Thus, we have
$$
d\tilde s^2=du^2+\sin^2\alpha\, dv^2+4\sin^4(\alpha/2)dv^2=du^2+4\sin^2(\alpha/2)dv^2.
$$
Taking into account \eqref{curv} we can easily find the Gaussian curvature of the totally geodesic submanifold
$\xi(M^2)$, namely
$$
\tilde K=\frac14 K(K-2\cot(\alpha/2)K'_\alpha),
$$
where $K(\alpha)$ is the Gaussian curvature of $M^2$ given by relations \eqref{curv} and \eqref{alf}.
\end{remark}

The equations \eqref{alf} and \eqref{curv} completely determine the class of
Riemannian 2-dimensional manifolds admitting a totally geodesic unit vector field.

\begin{proposition} Let $M^2$ be a Riemannian manifold with a line element of the form
$$
ds^2=du^2+\sin^2\alpha(u) dv^2.
$$
Denote by $K$ the Gaussian curvature of $M^2$. Then $\ds K=\frac{d\alpha}{du} $ if and only if the function
$\alpha(u)$ satisfies
$$
\frac{d\alpha}{du} =1+\frac{m}{\cos\alpha} \quad (m=const).
$$
\end{proposition}
\begin{proof}
The sufficient part is already proved. Suppose now that
$$
\frac{d\alpha}{du}=K (\ne0).
$$
Then we have
$$
\alpha'=K=-\frac{\partial_{uu}( \sin\alpha)}{\sin\alpha}=(\alpha')^2-\cot\alpha\, \alpha''.
$$
Therefore, $\ds \alpha''=-\alpha'(1-\alpha')\tan\alpha, $ \quad or \quad $ \ds
\frac{\alpha''}{\alpha'-1}=\alpha'\tan\alpha, $ \quad or
$$
(\ln|\alpha'-1|)'=-(\ln|\cos\alpha|)'.
$$
 Evidently, now
$\ds |\alpha'-1|=\frac{|m|}{|\cos\alpha|} $ where $m=const$ is a constant of integration. Finally,\ \ $\ds
\frac{d\alpha}{du}=1+\frac{m}{\cos\alpha}. $
\end{proof}

\begin{corollary}
Let $M^2$  be a Riemannian manifold of constant curvature $c\ne0$. Then $M^2$ admits a totally geodesic unit
vector field if and only if $c=1$. This vector field is parallel along meridians and moves along parallels with
unit angle speed.
\end{corollary}

\begin{proof} If $K=c=const$, then \eqref{curv} can be satisfied if and only if $c=1, \,a=-1$.
\end{proof}

The equation \eqref{curv} implies an elementary non-existence result.
\begin{corollary} Let $M^2$ be a Riemannian manifold with Gaussian curvature $K$. Then $M^2$ does not admit a
totally geodesic unit vector field $\xi$ with angle speed $a$ if $|K-1|<|a+1|$.
\end{corollary}
\begin{proof}
Indeed, one can easily see that $\ds \cos\alpha=\frac{a+1}{1-K}. $ If $|a+1|>|K-1|$,
then we come to a contradiction.
\end{proof}

\section{Integral trajectories of the totally geodesic vector field}

The integral trajectories of the totally geodesic vector field $\xi$ can be found easily as follows. Let
$\gamma=\{u(s),v(s)\}$ be an integral trajectory. Since
$$
\xi=\cos\omega\,e_0+\sin\omega\,e_1=\cos\omega\,\partial_u+\frac{\sin\omega}{\sin\alpha}\,\partial_v ,
$$
we can set
$$
\frac{du}{ds}=\cos\omega,\quad \frac{dv}{ds}=\frac{\sin\omega}{\sin\alpha}
$$
and then
$$
\frac{du}{dv}=\cot\omega \sin\alpha.
$$
Since $\alpha=\alpha(u)$ and $\omega=av+\omega_0$, we come to the equation with
separable variables
$$
\frac{du}{\sin\alpha}=\cot\omega\, dv.
$$
Using \eqref{alf}, we can find
$$\frac{du}{d\alpha}=\frac{\cos\alpha}{-a-1+\cos\alpha}
$$
and make a parameter change in the left hand side of the equation above. Then we come to the equation
$$
\frac{\cos\alpha\,d\alpha}{\sin\alpha(-a-1+\cos\alpha)}=\cot\omega\,dv.$$ Taking primitives, we have
$$
\begin{array}{ll}
\tan(\alpha/2)\sin(av+\omega_0)=c\,(a+(a+2)\tan^2(\alpha/2))^{\frac{a+1}{a+2}} &\quad \mbox{for \ }a\ne0,-2, \\[2ex]
\tan(\alpha/2)\sin(-2v+\omega_0)=c\,e^{\frac12\tan^2(\alpha/2)} &\quad \mbox{for \ } a=-2,\\[2ex]
\frac12\tan\omega_0\left(\frac{1}{1-\cos\alpha}+\ln\big|\tan(\alpha/2)\big|\right)=v-c &\quad \mbox{for \ }a=0.
\end{array}
$$
Taking into account \eqref{lambda}, we remark that $ \tan(\alpha/2)\sin\omega=k$ and
$\tan^2(\alpha/2)=k^2+\varkappa^2$. Therefore, we have an intrinsic equation on the integral curves of the
totally geodesic vector field
$$
\begin{array}{ll}
k=c\,\big[a+(a+2)(k^2+\varkappa^2)\big]^{\frac{a+1}{a+2}} &\quad \mbox{for \ }a\ne0,-2, \\[2ex]
k=c\,e^{\frac12(k^2+\varkappa^2)} &\quad \mbox{for \ } a=-2,\\[2ex]
k=\sin\omega_0\,\exp\big[2\cot\omega_0(v-c)-\frac12\frac{1+k^2+\varkappa^2}{k^2+\varkappa^2}\big] &\quad
\mbox{for \ } a=0,
\end{array}
$$
where $c$ is a constant of integration.

Moreover, in any case
\begin{multline*}
$$
\ds \xi(k)=\cos\omega\,\partial_u
[\tan(\alpha/2)\sin\omega]+\frac{\sin\omega}{\sin\alpha}\partial_v[\tan(\alpha/2)\sin\omega]=\\
\frac{\cos\omega\sin\omega\,\alpha'_u}{2\cos^2(\alpha/2)}+\frac{a\,\sin\omega\cos\omega\tan(\alpha/2)}{\sin\alpha}=
\frac{\cos\omega\sin\omega}{2\cos^2(\alpha/2)}(\alpha'_u+a).
$$
\end{multline*}
The equation \eqref{alf} yields
$$
\xi(k)=\frac{(a+1)\cos\omega\sin\omega}{2\cos^2(\alpha/2)}\left(1-\frac{1}{\cos\alpha}\right).
$$

Thus, if $a=-1$, then the integral trajectories of the field $\xi$ form a family of circles.  The metric of
$M^2$ is
$$
ds^2=du^2+\sin^2u\,dv^2
$$
 and we are dealing with the unit sphere parameterized by
$$
r=\big\{\sin u\cos v, \sin u\sin v, \cos u\big\}.
$$
These circles satisfy
\begin{equation}\label{lines}
\tan(u/2)\sin v=c.
\end{equation}
Let $(\rho,\varphi)$ be polar coordinates in a Cartesian plane which passes through the center of the sphere
such that $(0,0,1)$ is the \textit{north} pole on the sphere. Then the parameters $(\rho,\varphi)$ and $(u,v)$
are connected via stereographic projection from the \textit{south }pole as
$$
\begin{array}{l}
\rho=\tan(u/2),\\
\varphi=v.
\end{array}
$$
Therefore, the equation \eqref{lines} defines a family of parallel straight lines on the Cartesian plane.
\textit{The family of integral curves of a totally geodesic vector field on the unit sphere can be obtained as
inverse images under stereogrphic projection of this family.}

An explicit equation of this family is
$$
r(v)=\left\{ \frac{2c\sin v \cos v}{c^2+\sin^2v}, \frac{2c\sin^2v
}{c^2+\sin^2v},-\frac{c^2-\sin^2v}{c^2+\sin^2v} \right\}
$$
where  $c$ is the geodesic curvature of the corresponding circle. All of these circles pass through the south
pole $(0,0,-1)$ when $v=0,\pi$. We can find this by using the expression $\tan(u/2)=c/\sin v$ and trigonometric
expressions for $\sin u$ and $\cos u$ via $\tan(u/2)$.

The unit sphere is not the unique surface that realizes the metric \eqref{met1}.  Let $(x,y,z)$ be standard
Cartesian coordinates in $E^3$. We can find an isometric immersion of the metric \eqref{met1} into $E^3$ in a
class of  a surfaces of revolution. To do this, set
$$
\begin{array}{l}
\ds x(\alpha)=\sin\alpha,\\[1ex]
\ds (x'_\alpha)^2+(z'_\alpha)^2=\left(\frac{\cos\alpha}{a+1-\cos\alpha}\right)^2
\end{array}
$$
and we easily find
$$
\begin{array}{l}
\ds x(\alpha)=\sin\alpha,\\[2ex]
\ds z(\alpha)=\int_{\alpha_0}^\alpha\frac{\cos t }{a+1-\cos t }\sqrt{1-(a+1-\cos t
)^2}\ d\,t,
\end{array}
$$
where the interval of integration is limited by the restrictions
$$
\left\{
\begin{array}{l}
 1+a<\cos\alpha<2+a,\\
 -2<a<-1,
\end{array}
 \right.
 \quad \mbox{or}\quad
 \left\{
 \begin{array}{l}
 a<\cos\alpha<1+a,\\
-1<a<0.
\end{array}
\right.
$$
The restrictions mean that if $|a+1|\geq1$, then the metric \eqref{met1} does not
admit an isometric immersion into $E^3$ in a class of surfaces of revolution.

 \vspace{1cm}

\noindent
Department of Geometry,\\
Faculty of Mechanics and Mathematics,\\
Kharkiv National University,\\
Svobody Sq. 4,\\
 61077, Kharkiv,\\
Ukraine.\\
e-mail: yamp@univer.kharkov.ua


\begin{thebibliography}{30}

\bibitem{Am} Aminov Yu. The geometry of vector fields. Gordon\&Breach Publ., 2000.

\bibitem{BX-V1}
Boeckx~E., Vanhecke~L. \textit{Harmonic and minimal radial vector fields.} Acta Math. Hungar. 90 (2001),
317-331.

\bibitem{BX-V2}
Boeckx~E., Vanhecke~L. \textit{Harmonic and minimal vector fields on tangent and unit tangent bundles.}
Differential Geom. Appl. 13 (2000), 77-93.

\bibitem{Bx-V4}
Boeckx~E., Vanhecke~L. \textit{Isoparametric functions and harmonic and minimal unit vector fields}, Contemp.
Math. 288 (2001), 20--31.


\bibitem{B-Y} Borisenko A., Yampolsky A.\textit{ Riemannian geometry of bundles.}
Uspehi Mat. Nauk, 26/6 (1991), 51-95; Engl. transl.: Russian Math. Surveys, 46/6 (1991), 55--106.


\bibitem{B-Ch-N} Brito F., Chacon P., Naveira A.\textit{ On the volume of vector fields on spaces of constant sectional
 curvature}, Preprint, 2001.

\bibitem{GM-LF}
Gil-Medrano~O.,Llinares-Fuster~E. \textit{Minimal unit vector fields}, T\^ohoku Math. J. 54 (2002), 71 -- 84.

\bibitem{GM-LF2}Gil-Medrano~O.,Llinares-Fuster~E. \textit{Second variation of volume and energy of vector fields.
Stability of Hopf vector fields,} Math. Ann. 320 (2001), 531-545.

\bibitem{GM}Gil-Medrano~O. \textit{Relationship between volume and energy of unit vector fields}, Diff. Geom. Appl. 15
(2001), 137-152.

\bibitem{GM-GD-Vh} Gil-Medrano~O., Gonz\'alez-D\'avila J.C., Vanhecke L. \textit{Harmonic and minimal invariant unit
vector fields on homogeneous Riemannian manifolds}, Houston J. Math. 27 (2001), 377-409.

\bibitem{G-Z}
Gluck~H., Ziller~W. \textit{On the volume of a unit vector field on the three-sphere}. Comm. Math. Helv. 61
(1986), 177-192.

\bibitem{GD-V1}
Gonz\'alez-D\'avila J.C., Vanhecke L. \textit{Examples of minimal unit vector fields.} Ann Global Anal. Geom. 18
(2000), 385-404.

\bibitem{GD-V2}
Gonz\'alez-D\'avila J.C., Vanhecke L. \textit{Minimal and harmonic characteristic vector fields on
three-dimensional contact mentic manifolds}. J. Geom. 72 (2001), 65-76.

\bibitem{Ped}
Pedersen~S.~L. \textit{Volumes of vector fields on spheres}, Trans. Amer. Math. Soc. 336 (1993), 69-78.

\bibitem{S}Sasaki~S. \textit{On the differential geometry of tangent bundles of Riemannian manifolds.}
T\^ohoku Math. J. 10 (1958), 338-354.

\bibitem{TS-V1}
Tsukada K.,Vanhecke~L. \textit{Minimality and harmonicity for Hopf vector fields}, Illinois J. Geom. 45 (2001),
441 -- 451.

\bibitem{TS-V2}
Tsukada K.,Vanhecke~L. \textit{Invariant minimal unit vector fields on Lie groups,} Period. Math. Hungar. 40
(2000), 123-133.

\bibitem{TS-V3} Tsukada K.,Vanhecke~L. \textit{Minimal and harmonic vector fields on
$G_2(C^{m+2})$ and its dual space}. Monatsh. Math. 130 (2000), 143-154.

\bibitem{Ym1} Yampolsky A. \textit{On the mean curvature  of a unit vector field,} Math. Publ.
Debrecen, 60/1-2 (2002), 131 -- 155.

\bibitem{Ym2} Yampolsky A. \textit{On the intrinsic geometry of a unit vector field,}
Comment. Math. Univ. Carol. 43/2 (2002), 299-317.

\bibitem{Ym3} Yampolsky A. \textit{A totally geodesic property of Hopf vector fields,} Acta Math. Hungar. 101/1-2
(2003), 93-112.

\bibitem{Ym4} Yampolsky A. \textit{On extrinsic geometry of unit normal vector field of Riemannian hyperfoliation,}
 Math. Publ. Debrecen {63}/4 (2003), 555 -- 567.


\end{thebibliography}
\end{document}